\documentclass[12pt,a4paper]{article}
\usepackage{amstext}
\usepackage{fancyhdr}
\usepackage{amsfonts,graphicx,bezier, amssymb}

\usepackage{caption}
\usepackage{etoolbox}
\usepackage{pgfcore}
\usepackage{caption}
\usepackage{graphics,graphicx, bezier, float, color}
\usepackage{amsmath,amsthm,amssymb,amsfonts,amscd,enumerate, array,latexsym,epsfig,psfrag}
 \newtheorem{thm}{Theorem}[section]
  \newtheorem{defn}[thm]{Definition}
\let\olddefinition\defn
\renewcommand{\defn}{\olddefinition\normalfont}
 \newtheorem{prop}[thm]{Proposition}
  \newtheorem{lem}[thm]{Lemma}
  \newtheorem{cor}[thm]{Corollary}
  \newtheorem{exam}[thm]{Example}
\newenvironment{prf}{\noindent{\bf{Proof:}}~~}{\hfill\rule{1ex}{1ex}\vskip1.5ex}
\newcommand{\Z}{\mathbb Z}

\newcommand{\Ra}{\Rightarrow}

\newcommand{\beqa}{\begin{eqnarray}}
\newcommand{\enqa}{\end{eqnarray}}
\newcommand{\beq}{\begin{eqnarray*}}
\newcommand{\enq}{\end{eqnarray*}}

\begin{document}
\begin{center}
{\bf\Large Generalised reduced modules}
\vspace*{0.5cm}

%
%

\end{center}

\vspace*{0.5cm}
\begin{center}
Annet Kyomuhangi\footnote{Department of Mathematics, Busitema University\\
 P.O BOX 236, Tororo Uganda}$^{,}$
\footnote{Department of Mathematics, Makerere University\\
 P.O BOX 7062, Kampala Uganda} and David Ssevviiri$^{2}$$^{,}$
\footnote{Corresponding author}\\
E-mail: annet.kyomuhangi@gmail.com and david.ssevviiri@mak.ac.ug

 \end{center}

\vspace*{0.5cm}
\begin{abstract}
Let $R$ be a commutative unital ring, $a\in R$ and $t$ a positive integer.
$a^{t}$-reduced $R$-modules and universally $a^{t}$-reduced $R$-modules are defined and their properties given. Known (resp. new) results about reduced $R$-modules are retrieved (resp. obtained) by taking $t=1$ and results about reduced rings are deduced. 
   \end{abstract}


\vspace*{0.4cm}

{\bf MSC 2010} Mathematics Subject Classification: 16N60, 16S90, 13D07, 16E50, 16D80

\section{Introduction}
  \paragraph\noindent
 Throughout this paper, $R$ denotes a commutative ring with unity. All modules are left $R$-modules and the category of $R$-modules is denoted by $R$-Mod.
  A ring is \emph{reduced} if it has no non-zero nilpotent elements. Lee and Zhou in \cite{Leezhou} generalised this notion to modules by defining reduced modules.  An $R$-module $M$ is \emph{reduced} if for all $a\in R$ and $m\in M,$ $a^{2}m=0$ implies that $am=0.$ Reduced modules were further studied by \cite{Agay,Ann,Rege,DS} among others. In \cite{Ann}, an equivalent formulation of reduced modules was obtained in terms of the $a$-torsion functor $\Gamma_{a}(-),$ linking reduced modules to local cohomology, a very important topic in both commutative algebra and algebraic geometry. The endo-functor $\Gamma_{a}(-)$ on the category of $R$-modules associates to each $R$-module $M$ a submodule $\Gamma_{a}(M)$ given by $\Gamma_{a}(M):= \left \{m\in M ~| ~a^{k}m= 0 ~\text{for some } k\in \mathbb{Z}^{+} \right \}.$
   Note that for a ring $R$, $a\in R$ and $I$ the ideal of $R$ generated by $a,$  $\Gamma_{I}(M)= \Gamma_{a}(M).$ For properties of the $a$-torsion functor, see \cite{Brodmann, 24, Rohrer, Peter}.
 An $R$-module $M$ is \emph{$a$-reduced} if for $a\in R,$ $a\Gamma_{a}(M)=0,$ see \cite[Proposition 2.2]{Ann}.  An $R$-module $M$ is reduced if it is $a$-reduced for all $a\in R.$
 It was shown in \cite[Corollary 4.1]{Ann} that if $\mathcal{N}(R)$ is the nilradical of a ring $R,$ then $ \mathcal{N}(R)= \bigcup\limits_{a\in R}a\Gamma_{a}(R),$ stratifying the nilradical of $R.$ The functor $a\Gamma_{a}(-)$ on the category $R$-Mod was called the \emph{locally nilradical}.
\paragraph\noindent
 In this paper, a generalisation of the $a$-reduced (resp. reduced) $R$-modules is given by defining $a^{t}$-reduced (resp. universally $a^{t}$-reduced) $R$-modules for some positive integer $t$ and $a\in R.$  By taking $t=1,$ results that were obtained in \cite{Ann} are retrieved and new results about $a$-reduced (resp. reduced) $R$-modules are obtained. If $t$ is a positive integer and $a\in R,$ then an $R$-module $M$ is {\emph{$a^{t}$-reduced}} if for every positive integer $k\geq t$ and $m\in M,$ $a^{k}m=0$ implies that $a^{t}m=0.$ This is equivalent to saying that $a^{t}\Gamma_{a}(M)=0,$ (Proposition \ref{ak}). An $R$-module is \emph{universally $a^{t}$-reduced (written as $\varepsilon^{t}$-reduced for brevity)} if it is $a^{t}$-reduced for all $a\in R.$ We call the functor $a^{t}\Gamma_{a}: R\text{-Mod} \to R\text{-Mod}$ which associates to every $R$-module $M$ a submodule  $a^{t}\Gamma_{a}(M)$ of $M$ given by \{$a^{t}m~|~m\in \Gamma_{a}(M)$\}, the \emph{generalised locally nilradical}. In Section \ref{k}, properties of $a^{t}$-reduced (resp. universally $a^{t}$-reduced) $R$-modules have been given.

 \paragraph\noindent
 Regular rings were used in \cite[Theorem 2.16]{Rege} to characterise reduced modules. A ring is regular if and only if every $R$-module is reduced.
 Let $t$ be a positive integer. $t$-regular rings were defined in \cite{McCoy} as a generalisation of regular rings. In Section \ref{regular}, we show that every $R$-module is $\varepsilon^{t}$-reduced if and only if $R$ is a $t$-regular ring and all its ideals of the form $(0:b)$ for some $0\neq b\in R$ are semiprime, (Theorem \ref{reg}). It is deduced that a Noetherian ring is $t$-regular if and only if it is reduced, (Corollary \ref{noetred}).
   Section \ref{pro} is devoted to studying properties of the generalised locally nilradical $a^{t}\Gamma_{a}(-)$.  $a^{t}\Gamma_{a}(-)$ is a radical, (Proposition \ref{nilKrad}) which is not left exact. It is proved that for any ring $R,$  and $a\in R,$ $a^{t}\Gamma_{a}(R[x])=a^{t}\Gamma_{a}(R)[x],$ (Proposition \ref{poly}). As such, a polynomial ring $R[x]$ is $a^{t}$-reduced if and only if so is the ring $R,$ (Corollary \ref{pol}).
  \paragraph\noindent
  In Section \ref{Ser}, we show that $a^{t}$-reduced modules simplify the computation of local cohomology, i.e., if $R$ is a Noetherian ring, $a\in R,$ $I$ an ideal of $R$ generated by $a,$  and $M$ an $a^{t}$-reduced $R$-module, then the $i$-th local cohomology, $H_{I}^i(M)$ of $M$  is given by $H_{I}^i(M)\cong\text{Ext}_{R}^i(R/I^{t}, M),$ (Proposition \ref{tz}). 

  \paragraph\noindent
  Throughout this paper, $t$ is a fixed positive integer.
  Given an $R$-module $M,$ we write $N\leq M$ (resp. $I\lhd R$) to mean $N$ is a submodule of $M$ (resp. $I$ is an ideal of $R$).  For $a\in R$ and $k\in \Z^{+},$ $(0:_{M}{a}^{k})$ denotes a submodule of $M$ given by $\left\{m\in M~|~a^{k}m=0\right\}$ whereas, for any submodule $N$ of an $R$-module $M$ and $m \in M\setminus N,$ $(N: m):=\{r\in R~|~rm\in N\}.$

 \section{Properties of ${\varepsilon}^{t}$-reduced modules}\label{k}
 \paragraph\noindent
 Let $a\in R.$ An $R$-module $M$ is \emph{$a$-reduced} \cite[Definition 2.1]{Ann} if for all $m\in M,$ $a^{2}m=0$ implies that $am=0.$ An $R$-module is \emph{reduced} if it is $a$-reduced for all $a\in R.$  

 \begin{defn}\label{k-red}
 Let $a\in R$ and $t$ be a positive integer.
  An $R$-module $M$ is \emph{$a^{t}$-reduced} if for all $m\in M$ and all positive integers $k\geq t,$ $a^{k}m=0$ implies that $a^{t}m=0.$  Moreover, observe that for $t=1,$ we get exactly the notion of $a$-reduced module.
 \end{defn}
 \begin{defn}\label{t}

 An $R$-module $M$ is \emph{universally $a^{t}$-reduced} (written as \emph{${\varepsilon}^{t}$-reduced} for brevity) if for all $a\in R,$ $m\in M$ and all positive integers $k\geq t,$ $a^{k}m=0$ implies that $a^{t}m=0.$
 \end{defn}
   \paragraph\noindent Thus, an $R$-module $M$ is \emph{${\varepsilon}^{t}$-reduced} if and only if it is \emph{$a^{t}$-reduced} for all $a\in R.$ A ring $R$ is \emph{$a^{t}$-reduced} (resp. \emph{${\varepsilon}^{t}$-reduced}) if it is $a^{t}$-reduced (resp. ${\varepsilon}^{t}$-reduced) as an $R$-module. A submodule $N$ of an $R$-module $M$ is an \emph{$a^{t}$-reduced  (resp. ${\varepsilon}^{t}$-reduced) submodule} if it is $a^{t}$-reduced (resp. ${\varepsilon}^{t}$-reduced) as an $R$-module.

 \begin{prop}\label{ak}
 Let $a\in R,$ $I$ be the ideal of $R$ generated by $a$ and $M$ an $R$-module. The following statements are equivalent:
 \begin{enumerate}
 \item $M$ is $a^{t}$-reduced;
 \item $a^{t}\Gamma_{a}(M)=0;$
  \item $(0:_{M}{a}^{k})=(0:_{M}{a}^{t})$ for all positive integers $k\geq t;$
\item $\text{ Hom}_{R}(R/I^{k},M) \cong \text{ Hom}_{R}(R/I^{t},M)$ for all positive integers $k\geq t;$
\item $\Gamma_{a}(M)\cong \text{ Hom}_{R}(R/I^{t},M);$

\item $0\to \Gamma_{a}(M) \to M \to a^{t}M \to 0 $ is a short exact sequence.
 \end{enumerate}
 \end{prop}

 \begin{prf}
 \begin{enumerate}
 \item[$1\Ra 2$] Let $m \in a^{t}\Gamma_{a}(M).$ $m=a^{t}n$ for some $n\in \Gamma_{a}(M).$ So, there exists $k\in \Z^{+}$ such that $a^{k}n=0.$ If $k<t,$ then $a^{t}n=0$ and we are done. Suppose $ k\geq t.$ By $1,$ $a^{t}n=0.$ So, $a^{t}\Gamma_{a}(M)=0.$

 \item[$2\Ra 3$] For $k\geq t,~ (0:_{M}{a}^{t})\subseteq(0:_{M}{a}^{k})$. If $m\in (0:_{M}{a}^{k}),$ then $a^{k}m=0.$ This implies that $m\in \Gamma_{a}(M).$ By $2,$ $a^{t}m=0.$ This shows that $m\in (0:_{M}{a}^{t})$ which completes the proof.
 \item[$3\Ra 4$] Since $I$ is generated by $a,$ $\text{Hom}_{R}(R/I^{k},M)\cong(0:_{M}{a}^{k}).$ It follows from $3$ that $\text{Hom}_{R}(R/I^{k},M)\cong \text{Hom}_{R}(R/I^{t},M)$ for all positive integer $k\geq t.$
 \item[$4\Ra 5$] $\Gamma_{a}(M) \cong \underset{k}\varinjlim\text { Hom}_{R}(R/I^{k},M),$ see \cite[page $6$]{Brodmann}. By $4,$ $\Gamma_{a}(M)\cong \text {Hom}_{R}(R/I^{t},M).$

  \item[$5\Ra 6$] The $R$-homomorphism $f: M \to a^{t}M$ defined by $m\mapsto a^{t}m $ is an epimorphism. Ker$f=(0:_{M}{a}^{t})\cong \text{Hom}_{R}(R/I^{t},M).$ By $5,$ Ker$f= \Gamma_{a}(M).$ So, $0\to \Gamma_{a}(M) \to M \to a^{t}M \to 0 $ is a short exact sequence.
 \item[$6\Ra 1$] Suppose $a\in R$ and $m\in M$ such that for some positive integer $k\geq t,$ $a^{k}m=0.$ This implies that $m\in \Gamma_{a}(M).$ By $6,$ $\Gamma_{a}(M)= (0:_{M}{a}^{t}).$ So, $a^{t}m=0.$ By Definition \ref{k-red}, $M$ is $a^{t}$-reduced.

 \end{enumerate}
 \end{prf}

 \paragraph\noindent
It is well known that the $a$-torsion functor $\Gamma_{a}(-)$ is idempotent and $R$-linear. It follows that $$a\Gamma_{a}(a\Gamma_{a}(M))=a^{2}\Gamma_{a}(\Gamma_{a}(M))=a^{2}\Gamma_{a}(M).$$ Therefore, composing the functor $a\Gamma_{a}(-)$ $t$-times is nothing but the functor $a^{t}\Gamma_{a}(-).$ This provides another view point of the generalised locally nilradical $a^{t}\Gamma_{a}(-).$

\begin{cor}\label{red.}
 Let $a\in R,$ $I$ be the ideal of $R$ generated by $a,$ $M$ an $R$-module and $t\in \Z^{+}.$ The following statements are equivalent:
 \begin{enumerate}
\item $M$ is ${\varepsilon}^{t}$-reduced;
\item For all $a\in R,$ $a^{t}\Gamma_{a}(M)=0;$
\item $(0:_{M}{a}^{k})=(0:_{M}{a}^{t})$ for all positive integers $k\geq t$ and all $a\in R;$
\item $\text{Hom}_{R}(R/I^{k},M) \cong \text{Hom}_{R}(R/I^{t},M)$ for all ideals $I$ and all positive integers $k\geq t;$
\item $\Gamma_{a}(M)\cong \text{Hom}_{R}(R/I^{t},M)$ for all $a\in R;$
\item For all $a\in R,$ $0\to \Gamma_{a}(M) \to M \to a^{t}M \to 0 $ is a short exact sequence.
 \end{enumerate}
 \end{cor}
 \begin{prf}
 Follows from Proposition \ref{ak} and the fact that $M$ is ${\varepsilon}^{t}$-reduced if it is $a^{t}$-reduced for all $a\in R.$
 \end{prf}
 \paragraph\noindent
 Proposition \ref{ak} (resp. Corollary \ref{red.}) shows that, just like $a$-reduced (resp. reduced) modules, $a^{t}$-reduced (resp. ${\varepsilon}^{t}$-reduced) modules are categorical.
\begin{exam}
For $a\in R,$ $a^{t}$-reduced modules include:
\begin{enumerate}
  \item every $a$-torsion-free $R$-module, i.e., modules $M$ for which $\Gamma_{a}(M)=0,$
  \item a flat (resp. free, projective) module over an $a^{t}$-reduced ring $R.$
  \end{enumerate}
\end{exam}

\begin{prop}\label{inclu}
For $a\in R,$ we have\\
\begin{center}
   $
\begin{tabular}{cccc}
    reduced $R$-module & $\Rightarrow$ & $a$-reduced $R$-module\\
   \rotatebox[origin=c]{270}{$\Rightarrow$}& & \rotatebox[origin=c]{270}{$\Rightarrow$}\\
    ${\varepsilon}^{t}$-reduced $R$-module & $\Rightarrow$ & $a^{t}$-reduced $R$-module\\
\end{tabular}
.$
\end{center}
\end{prop}
\begin{prf}
 Let $a\in R$ and $M$ be an $R$-module. Since $a^{t}\Gamma_{a}(M)\subseteq a\Gamma_{a}(M),$  $a\Gamma_{a}(M)=0$ implies that $a^{t}\Gamma_{a}(M)=0.$ So, $a$-reduced (resp. reduced) implies  $a^{t}$-reduced (resp. ${\varepsilon}^{t}$-reduced). The other two implications follow immediately from the definitions.
\end{prf}
\paragraph\noindent
The converse of the implications in Proposition \ref{inclu} does not hold in general.
 Let $p$ be a prime number and $\Z_{p^{t}}$ be a group of integers modulo $p^{t}.$
Every $\Z$-module $\Z_{p^{t}}$ is $p^{t}$-reduced since $p^{t}\Gamma_{p}(\Z_{p^{t}})= 0.$ However, for $t>1,$  $p\Gamma_{p}(\Z_{p^{t}})\neq 0.$ 
  The $\Z$-module $\Z_{16}$ is $4^{2}$-reduced but not $2^{2}$-reduced and therefore not $\varepsilon^{t}$-reduced. Not every $a$-reduced module is reduced as the $\Z$-module $\Z_{4}$ is $4$-reduced but not reduced.
  If $R:=\Z_{8},$ $M:=\Z_{4}$ and $t:=2\in \Z^{+},$ then $\Z_{4}$ is $\varepsilon^{2}$-reduced but not reduced.

\begin{prop}
Let $S$ be a multiplicatively closed set of a ring $R.$ An $R$-module $M$ is ${\varepsilon}^{t}$-reduced if and only if the localisation module $S^{-1}(M)$ of $M$ over the ring $S^{-1}(R)$ is ${\varepsilon}^{t}$-reduced.
\end{prop}
\begin{prf}
Let $a \in R,$ $r,s\in S$ and $m\in M$ such that $\displaystyle{\frac{a}{r}\in S^{-1}(R)}$ and $\displaystyle{\frac{m}{s}\in S^{-1}(M)}.$
Suppose $M$ is $a^{t}$-reduced for all $a\in R$ and $\displaystyle{\bigg(\frac{a}{r}\bigg)^{k}\bigg(\frac{m}{s}\bigg)=0}$ for $k\geq t.$  So, $a^{k}m=0.$ Since $M$ is $a^{t}$-reduced for all $a\in R,$ $a^{t}m=0.$ Thus, $\displaystyle{ \frac{a^{t}m}{r^{t}s}=\bigg(\frac{a}{r}\bigg)^{t}\bigg(\frac{m}{s}\bigg)=0}$ and $S^{-1}(M)$ is ${\varepsilon}^{t}$-reduced.
The converse holds since $M$ can be embedded in $S^{-1}(M)$ and ${\varepsilon}^{t}$-reduced modules are closed under submodules.
\end{prf}
\begin{prop}
If for all $R$-modules $M$ and $N\leq M,$ the $R$-module $R/(N: m)$ for $m \in M\setminus N$ can be embedded in $M,$ then $a^{t}$-reduced $R$-modules are closed under epimorphic images.
\end{prop}
\begin{prf}
Let $M$ be an $a^{t}$-reduced $R$-module and $N$ a submodule of $M.$ By Proposition \ref{ak}, $a^{t}\Gamma_{a}(M)=0.$ By hypothesis, $R/(N: m)$ is embeddable in $M.$ As such, $a^{t}\Gamma_{a}\bigg(\displaystyle{\frac{R}{(N: m)}}\bigg)=0.$  This shows that $R/(N: m)$ is $a^{t}$-reduced as an $R$-module. So, if $a^{k}m\in N,$ then $a^{t}m\in N$ and $M/N$ is $a^{t}$-reduced.
\end{prf}

\section{Characterisation in terms of $t$-regular rings}\label{regular}
\paragraph\noindent
A ring $R$ is regular if and only if every $R$-module is reduced. In this section, we use $t$-regular rings to characterise ${\varepsilon}^{t}$-reduced modules.

\begin{defn}\cite[Page 175]{McCoy}\label{McCoy}
 A ring $R$ is \emph{$t$-regular} if for every $a\in R,$ there exists $b\in R$ such that $a^{t}=a^{2t}b,$ i.e., if every element $a^{t}$ of $R$ is regular.
\end{defn}
\paragraph\noindent
By \cite{Azumaya}, Definition \ref{McCoy} is equivalent to the statement: for every $a\in R,$ there exists $b\in R$ such that $a^{t}=a^{t+1}b.$

 \begin{lem}\label{**}
 $t$-regular rings are closed under quotients.
 \end{lem}
 \begin{prf}
 Suppose $R$ is a $t$-regular ring and $I$ is an ideal of $R.$  Let $\bar{a}\in R/I,$ where $\bar{a}= a+I$ for $a\in R.$ Thus, $\bar{a}^{t}= a^{t}+I.$ But $a^{t}=a^{2t}b.$ So, $\bar{a}^{t}= a^{2t}b+I= (a^{2t}+I)(b+I)=\bar{a}^{2t}\bar{b}.$
 \end{prf}
  \paragraph\noindent
  A regular ring is $t$-regular for any $t\in \Z^{+}.$  Special primary rings (rings whose elements are either nilpotent or units) defined in \cite{McCoy} are also $t$-regular. In general, a $t$-regular ring need not be regular. However, a $1$-regular ring is regular.
  \begin{prop}

  A domain is $t$-regular if and only if it is a field.
  \end{prop}
  \begin{prf}
  Suppose $R$ is $t$-regular and a domain. If $a\in R,$ then $a^{t}=a^{2t}b$ for some $b\in R.$ So, $a^{t}(1-a^{t}b)=0.$ Since $R$ is a domain, $a^{t}=0$ or $a^{t}b=1$ which implies that $a=0$ or $a$ is a unit. Thus, all non-zero elements of $R$ are units.
  For the converse, suppose that $a^{t}\in R$ and $R$ is a field. There exists $0\neq b\in R$ such that $a^{t}b=1.$ So, $a^{t}ba^{t}=a^{t}$ and $R$ is $t$-regular.  
  \end{prf}

  \paragraph\noindent
  A regular ring is reduced. For $t$-regular rings, we have Proposition \ref{***}.

    \begin{prop}\label{***}
   If $R$ is a $t$-regular ring such that every ideal of $R$ of the form $(0:b)$ for some $0\neq b\in R$ is semiprime, then $R$ is ${\varepsilon}^{t}$-reduced.
   \end{prop}
  \begin{prf}
  Suppose $a\in R$ such that $a^{k}r=0$ for $r\in R$ and a positive integer $k>t$ but $a^{t}r\neq 0.$ Since $R$ is $t$-regular, $a^{t}=a^{2t}b$ for some $0\neq b\in R$  and $a^{2t}br\neq 0.$ So, $a\notin (0:a^{2t-1}br).$ By hypothesis, $(0:a^{2t-1}br)$ is a semiprime ideal of $R.$  So, $a^{k}\notin (0:a^{2t-1}br)$ for all $k\in \Z^{+}$ and $a^{k}r\neq 0$ which is a contradiction.
  \end{prf}

  \begin{cor}\label{Noet}

  A Noetherian $t$-regular ring is ${\varepsilon}^{t}$-reduced.
  \end{cor}

  \begin{prf}
    For a Noetherian ring $R,$ a set $Z(R)$ of all zero-divisors of $R$ coincides with $\bigcup\limits_{P\in \text{Ass}(R)}P$ the union of all associated primes $P$ of $R,$ see \cite[Theorem 3.1]{Eis}. However, every prime ideal is semiprime. So, $Z(R)=\bigcup\limits_{0\neq b\in R }(0:b)$ consists of only semiprime ideals. The rest follows from Proposition \ref{***}.
   \end{prf}

   \begin{prop}\label{coi}
    $R$ is a regular ring if and only if it is $t$-regular and every ideal of $R$ is semiprime.
   \end{prop}
\begin{prf}
      Suppose $R$ is $t$-regular.  Let $J \lhd R$ such that $a\in J.$ It follows that  $a^{t}=(a^{t}b)a^{t}\in J^{2}.$ Since each ideal of $R$ is semiprime, $a \in J^{2}$ and $J=J^{2},$ i.e., each ideal of $R$ is idempotent. By \cite[Theorem 1.16]{Goodearl}, the proof is complete if we can show that every simple $R$-module is injective. This is exactly what was done in the proof $(d\Ra a)$ of \cite[Theorem 1.16]{Goodearl}.
      For completeness, we repeat that proof here. Suppose $I$ is a maximal ideal and $f: J \to R/I $ is a non-zero homomorphism. Since every ideal of $R$ is idempotent, $I\cap J=(I\cap J)^{2}\subseteq IJ\subseteq$ Ker$f \subset J.$ So, $J\nsubseteq I.$ Thus, there exists $i \in I$ and $j \in J$ such that $i+j=1$ and $a-aj=ai \in IJ\subseteq$ Ker$f.$ So, $f(a)=f(aj)=af(j)$ and $f$ extends to a map $R\to R/I.$  Conversely, if $R$ is regular, then by definition, it is $t$-regular. Moreover, by  \cite[Corollary 1.2]{Goodearl}, every ideal of $R$ is semiprime.

   \end{prf}

  \begin{prop}\label{****}
  Let $f:R \to S$ be a ring homomorphism and $M$ an $S$-module. For $r\in R$ and $m\in M,$ $M$ is an $R$-module via $rm:=f(r)m.$
  \begin{enumerate}
  \item If $_{S}M$ is ${\varepsilon}^{t}$-reduced, then so is $_{R}M;$
  \item if $f$ is surjective and $_{R}M$ is ${\varepsilon}^{t}$-reduced, then so is $_{S}M.$
  \end{enumerate}
  \end{prop}
  \begin{prf}
  \begin{enumerate}
  \item Suppose $m\in M$ such that for all $a\in R,$ $a^{k}m=0$ for $k\geq t.$ Then, $f(a^{k})m= (f(a))^{k}m=0.$ By hypothesis, $(f(a))^{t}m=0.$ So, $a^{t}m=0$ and $_{R}M$ is ${\varepsilon}^{t}$-reduced.
  \item Let $s^{k}m=0$ for $s\in S$ and $k\in \Z^{+}.$ By hypothesis, $s=f(a)$ for some $a\in R$ and  $s^{k}m=(f(a))^{k}m= 0.$  So, $a^{k}m=0.$ Since $_{R}M$ is ${\varepsilon}^{t}$-reduced, $a^{t}m=0$ and $(f(a))^{t}m= 0.$ So, $s^{t}m=0$ and $_{S}M$ is  ${\varepsilon}^{t}$-reduced.
  \end{enumerate}
  \end{prf}
\begin{prop}\label{*}
Let $a\in R$ and $M$ be an $R$-module. The following statements are equivalent:
\begin{enumerate}
\item $M$ is ${\varepsilon}^{t}$-reduced,
\item every cyclic submodule of $M$ is ${\varepsilon}^{t}$-reduced.
\end{enumerate}
\end{prop}
\begin{prf}
\begin{enumerate}
\item[$1\Ra 2$] Follows from the fact that ${\varepsilon}^{t}$-reduced modules are closed under submodules.
\item[$2\Ra 1$] Suppose every cyclic submodule of $M$ is ${\varepsilon}^{t}$-reduced and $m\in M$ such that for all $a\in R$ and some positive integer $k\geq t,$ $a^{k}m=0.$ So, $a^{k}rm=0$ for all $r\in R.$ Since $rm\in Rm$ and $Rm$ is ${\varepsilon}^{t}$-reduced, $a^{t}rm=0.$  As such, $a^{t}Rm=0$ and $a^{t}m=0.$
    \end{enumerate}
\end{prf}

\begin{thm}\label{reg}
Let every ideal of $R$ of the form $(0:b)$ for some $0 \neq b\in R$ be semiprime.
The following statements are equivalent:
\begin{enumerate}
\item every $R$-module is ${\varepsilon}^{t}$-reduced,
\item every cyclic $R$-module is ${\varepsilon}^{t}$-reduced,
\item $R$ is $t$-regular.
\end{enumerate}
\end{thm}
\begin{prf}
\begin{enumerate}
\item[$1\Leftrightarrow 2$] Follows from Proposition \ref{*}.
\item[$2 \Ra 3$] By $2,$ the cyclic $R$-module $R/a^{2t}R$ is $a^{t}$-reduced for all $a\in R.$ It follows that  $a^{2t}.\bar{1}=0$ and $\bar{a}^{t}=a^{t}.\bar{1}=\bar{0}$ where $\bar{a}=a+a^{2t}R$ for any $a\in R.$ So, $a^{t} \in a^{2t}R $ and $a^{t} = a^{2t}b $ for some $b\in R.$
\item [$3\Ra 2$] Let $M$ be a cyclic $R$-module. Then $M\cong R/I$ for some ideal $I$ of $R.$ Since $R$ is $t$-regular, by Lemma \ref{**}, $R/I$ is also $t$-regular. By Proposition \ref{***}, $R/I$ is an ${\varepsilon}^{t}$-reduced ring. By Proposition \ref{****}$(1),$ the $R$-module $M\cong R/I$ is ${\varepsilon}^{t}$-reduced.

\end{enumerate}
\end{prf}

\begin{prop}\label{noeth}
 Let $R$ be a ring such that for every $R$-module $M$ and $0\neq m\in M,$ $(0:m)$ is a semiprime ideal of $R.$ $M$ is reduced if and only if it is ${\varepsilon}^{t}$-reduced.
 \end{prop}
\begin{prf}
By Proposition \ref{inclu}, every reduced module is ${\varepsilon}^{t}$-reduced. Conversely, if
$M$ is ${\varepsilon}^{t}$-reduced such that for $m\in M,$ $a\in R$ and $k\in \Z^{+},$ $a^{k}m=0,$ then $a^{k} \in (0:m).$ By hypothesis, $a \in (0:m).$ So, $am=0.$

\end{prf}
\begin{cor}\label{noth}

A finitely generated module over a Noetherian ring is reduced if and only if it is ${\varepsilon}^{t}$-reduced.
\end{cor}

\begin{cor}\label{nore}
 A Noetherian ring is reduced if and only if it is ${\varepsilon}^{t}$-reduced.
\end{cor}
\begin{prf}
A Noetherian ring $R$ has ideals of the form  $(0:b)$ for some $0\neq b\in R$ semiprime. The rest follows from Corollary \ref{noth}.
\end{prf}
\begin{cor}\label{noetred}
A Noetherian ring is $t$-regular if and only if it is reduced.
\end{cor}
\begin{prf}
Follows from Corollaries \ref{Noet} and \ref{nore}.
\end{prf}
\paragraph\noindent We recall that an $R$-module $M$ is \emph{cogenerated} by a ring $R$ if it can be embedded in a direct product of copies of $R.$ An $R$-module $M$ is \emph{faithful} if for each $a\in R,$  $aM=0$ implies that $a=0.$
\begin{prop}\label{Aga}
For $a\in R,$ the following statements are equivalent:
\begin{enumerate}
\item $R$ is $a^{t}$-reduced (resp. $\varepsilon^{t}$-reduced),
\item every $R$-module cogenerated by $R$ is $a^{t}$-reduced (resp. $\varepsilon^{t}$-reduced),
\item every submodule of a free $R$-module is $a^{t}$-reduced (resp. $\varepsilon^{t}$-reduced),
\item there exists a faithful $a^{t}$-reduced (resp. $\varepsilon^{t}$-reduced) $R$-module.
\end{enumerate}
\end{prop}
\begin{prf}
We only give a proof for the $a^{t}$-reduced case since the $\varepsilon^{t}$-reduced case follows trivially.
\begin{enumerate}
\item[$1\Ra 2$] Suppose $M$ is an $R$-module cogenerated by $R.$  $M$ can be embedded in some $\prod\limits_{i}R_{i}$ where $R_{i}=R.$ So, $M$ is $a^{t}$-reduced since a family of $a^{t}$-reduced modules is closed under submodules and direct products.
  \item[$2\Ra 3$] Suppose $F$ is a free $R$-module, i.e, $F\cong R^{n}$ for $n\in \Z^{+}.$  Every submodule of $F$ is cogenerated by $R.$ So, it is $a^{t}$-reduced. 
   \item[$3\Ra 4$] $R$ as an $R$-module is 
        faithful and a submodule of a free module. By $3,$ $R$ is an $a^{t}$-reduced $R$-module.
    \item[$4\Ra 1$] Let $M$ be a faithful $a^{t}$-reduced $R$-module. If $a,r \in R$ and $k\geq t \in \Z^{+}$ such that $a^{k}r=0,$ then for all $m\in M,$ $a^{k}rm=0.$ By hypothesis, $a^{t}rm=0$ for all $m\in M.$ So, $a^{t}rM=0.$  Since $M$ is faithful, $a^{t}r=0.$
\end{enumerate}
\end{prf}
\paragraph\noindent
 By writing $t=1$ in Proposition \ref{Aga}, we retrieve \cite[Theorem 2.6]{Agay} for an identity endomorphism.

\begin{cor}\label{cogen}

An ${\varepsilon}^{t}$-reduced cogenerator ring is $t$-regular.
\end{cor}
\begin{prf}
 Let Cog$(R)$ be a class of $R$-modules cogenerated by $R.$ Since Cog$(R)=R$-Mod, every $R$-module is ${\varepsilon}^{t}$-reduced. It follows from the proof $1\Leftrightarrow 2 \Ra 3$ of Theorem \ref{reg} that $R$ is $t$-regular.
\end{prf}
\paragraph\noindent
Every regular ring is reduced. In Corollary \ref{rr}, we give a condition for the converse to hold.
\begin{cor}\label{rr}
A reduced cogenerator ring is regular.
\end{cor}
\begin{prf}
Follows from Corollary \ref{cogen} by taking $t=1.$
\end{prf}

%
 \section{The generalised locally nilradical}\label{pro}

 \paragraph\noindent
 Let $a\in R$ and $M$ be an $R$-module. In this section, we give properties of the functor $a^{t}\Gamma_{a}:R{\normalfont\text{-Mod}} \to R{\normalfont\text{-Mod}}$ given by $M\mapsto a^{t}\Gamma_{a}(M) $ which we call the \emph{generalised locally nilradical}.
   \begin{paragraph}\noindent
A functor $\gamma:R\text{-Mod} \to R\text{-Mod} $ is a {\emph{preradical}} if for every $R$-homomorphism $f:M \to N,~ f(\gamma(M))\subseteq \gamma(N).$ A preradical
$\gamma $ is a {\emph{radical}} if for all $M \in R$-Mod,  $\gamma(M/\gamma(M))=0.$
 If $\gamma(N)= N \cap \gamma(M)$ for every  $N\leq M \in R$-Mod, then the radical $\gamma$ is called {\emph{hereditary}} (equivalently, {\emph{left exact}} as a functor).
\end{paragraph}

     \begin{prop}\label{nilKrad}
	For $a\in R,$ the functor $$a^{t}\Gamma_{a}:R{\normalfont\text{-Mod}} \to R{\normalfont\text{-Mod}}$$
 $$M\mapsto a^{t}\Gamma_{a}(M) $$ is a radical.
	\end{prop}
\begin{prf}
For an $R$-homomorphism $f: M \to N$ and $x\in f(a^{t}\Gamma_{a}(M)),$ $x=a^{t}f(m)$ for some $m\in \Gamma_{a}(M).$
$m\in \Gamma_{a}(M)$ implies that there exists $k\in \Z^{+}$ such that $a^{k}m=0.$ Thus, $a^{k}f(m)= f(a^{k}m)=f(0)=0.$ So, $f(m)\in \Gamma_{a}(N) $ and $x=a^{t}f(m)\in a^{t}\Gamma_{a}(N).$ Hence, $f(a^{t}\Gamma_{a}(M))\subseteq a^{t}\Gamma_{a}(N)$ which proves that $a^{t}\Gamma_{a}(-)$ is a preradical. By \cite[Proposition 3.1]{Ann}, $a\Gamma_{a}(M)$ is a radical and $a^{t}\Gamma_{a}(M)\subseteq a\Gamma_{a}(M).$ By \cite[Section 1.1.E1]{Bican}, we have $a\Gamma_{a}(M/a^{t}\Gamma_{a}(M))=a\Gamma_{a}(M)/a^{t}\Gamma_{a}(M).$ By left multiplication of $a^{t-1},$ we obtain  $a^{t}\Gamma_{a}(M/a^{t}\Gamma_{a}(M))=a^{t}\Gamma_{a}(M)/a^{t}\Gamma_{a}(M)=0.$
\end{prf}

\paragraph\noindent The functor $a^{t}\Gamma_{a}(-)$ in general is not left exact. For if  $R:=\Z,$ $M:=\Z_{8},$ $N:=2\Z_{8},$ $a:= 2\in \Z$ and $t:=2,$ then
 $2^{2}\Gamma_{2}(M)= 4\Z_{8}$ and $2^{2}\Gamma_{2}(N)= 0 \neq 4\Z_{8} =N \cap 2^{2}\Gamma_{a}(M).$ 
\paragraph\noindent A submodule $N$ of an $R$-module $M$ is a \emph{characteristic} submodule if for all automorphisms $f$ of $M$, $f(N)\subseteq N.$

\begin{prop}\label{proj}
	Let $a\in R$ and $M$ be an $R$-module. The following statements hold:
	\begin{enumerate}
		\item $a^{t}\Gamma_{a}(R)$ is an ideal of $R;$
		\item for each $M \in R$-Mod, $a^{t}\Gamma_{a}(M)$ is a characteristic submodule of $M$ and $$a^{t}\Gamma_{a}(R)M \subseteq a^{t}\Gamma_{a}(M);$$
		\item if $M$ is projective, then $a^{t}\Gamma_{a}(M)= a^{t}\Gamma_{a}(R)M.$
	\end{enumerate}

\end{prop}
\begin{prf}
Since $a^{t}\Gamma_{a}(-)$ is a preradical, the proof follows from \cite[Proposition 1.1.3]{Bican}.

\end{prf}

 \begin{prop}\label{factor}
	Let $a\in R,$ $M$ be an $R$-module and $N\leq M.$
	
		 $$a^{t}\Gamma_{a}(N) \subseteq N\cap a^{t}\Gamma_{a}(M)$$  and $$\frac{a^{t}\Gamma_{a}(M) + N}{N} \subseteq a^{t}\Gamma_{a}(M/N).$$
		
\end{prop}
\begin{prf}
	Follows from \cite[Proposition 1.1.1]{Bican} since $a^{t}\Gamma_{a}(-)$ is a (pre)radical.
\end{prf}
\begin{prop}\label{product}
	Let $I$ be an indexing set and $\{M_{i}\}_{i\in I}$ be a family of $R$-modules.
	$$ a^{t}\Gamma_{a}\left(\bigoplus\limits_{i\in I} M_{i}\right)=\bigoplus\limits_{i\in I}a^{t}\Gamma_{a}(M_{i}) $$ and $$ a^{t}\Gamma_{a}\left(\prod\limits_{i\in I} M_{i}\right) \subseteq \prod\limits_{i\in I}a^{t}\Gamma_{a}(M_{i}).$$
\end{prop}
\begin{prf}
	It follows from \cite[Proposition 1.1.2]{Bican}.
\end{prf}
\paragraph\noindent
If $M_{i}$ is $a^{t}$-reduced for each $i \in I,$ then $a^{t}\Gamma_{a}(M_{i})=0$ for each $i \in I.$ By Proposition \ref{product}, $a^{t}\Gamma_{a}\left(\prod\limits_{i\in I} M_{i}\right)=0.$  So, $\prod\limits_{i\in I} M_{i}$ is $a^{t}$-reduced which shows that a family of $a^{t}$-reduced modules is closed under direct products.

\begin{prop}\label{poly}
	For $a\in R,$  $$a^{t}\Gamma_{a}(R)[x]= a^{t}\Gamma_{a}(R[x]).$$
\end{prop}
\begin{prf}
Follows from the fact that the $a$-torsion functor $\Gamma_{a}(-)$ is $R$-linear and $a\Gamma_{a}(R)[x]= a\Gamma_{a}(R[x]),$ see \cite[Theorem 3.2]{Ann}.
\end{prf}

\begin{cor}\label{pol}
A polynomial ring $R[x]$ is $a^{t}$-reduced (resp. ${\varepsilon}^{t}$-reduced) if and only if $R$ is an $a^{t}$-reduced (resp. ${\varepsilon}^{t}$-reduced) ring.
\end{cor}
\begin{prf}
  By Proposition \ref{poly}, $a^{t}\Gamma_{a}(R)=0$ if and only if $a^{t}\Gamma_{a}(R[x])=0.$ The rest follows from Proposition \ref{ak}.
\end{prf}

\section{Application}\label{Ser}
\paragraph\noindent
Let $I\lhd R.$  Over Noetherian rings, the $I$-torsion functor $\Gamma_{I}(-)$ on $R$-Mod is a left exact radical. Its right derived functor which is known as the $i$-th local cohomology functor with respect to $I$ is defined by $H_{I}^{i}(-)\cong\underset{k}\varinjlim\text { Ext}_{R}^{i}(R/I^{k},-).$ For more information about local cohomology, see \cite{Brodmann, 24, Peter}.

\begin{prop}\label{tz}
If $R$ is a Noetherian ring, $a\in R,$  $I$ the ideal of $R$ generated by $a$ and $M$ an $R$-module, then the following statements hold:
\begin{enumerate}
\item if $M$ is an $a^{t}$-reduced $R$-module, then the $i$-th local cohomology module $H_{I}^i(M)$ is given by
\begin{equation*}
H_{I}^i(M)\cong\text{Ext}_{R}^i(R/I^{t}, M);
\end{equation*}

\item if $M$ is an $a^{t}$-reduced $R$-module and $R/I^{t}$ is a projective $R$-module, then for all $i\geq 1,$
$$H_{I}^{i}(M)=0.$$

\end{enumerate}

\end{prop}

 \begin{prf}
 \begin{enumerate}
\item  If $M$ is an $a^{t}$-reduced $R$-module, then by Proposition \ref{ak},
 $\Gamma_{I}(M)\cong\text{ Hom}_{R}(R/{I^{t}}, M)$. Thus, $H_{I}^i(M)\cong\text{Ext}_{R}^i(R/I^{t}, M)$.
  \item If $M$ is an $a^{t}$-reduced $R$-module, then by $1,$ $H_{I}^i(M)\cong\text{Ext}_{R}^i(R/I^{t}, M).$ Since $R/I^{t}$ is  projective, it follows by general theory that the cohomology module $H_{I}^i(M)$ vanishes for all $i\geq 1$.

 \end{enumerate}

 \end{prf}

~~

 {\bf{Acknowledgement}}

 \paragraph\noindent
 We acknowledge support from Sida bilateral programme (2015--2020) with Makerere University; Project 316: Capacity building in Mathematics and its applications and the Abram Gannibal project: Collaborative research in applied algebra and geometry in Africa.

\end{document}